\documentclass[10pt,a4paper]{article}
\usepackage{amsthm,amsfonts,amsmath,amscd,amssymb}
\usepackage[T1]{fontenc}

\usepackage[dvipsnames]{xcolor}
\usepackage{hyperref}
\usepackage{enumerate}
\usepackage{enumitem}
\usepackage{tikz}
\usepackage{tikz-cd}
\usepackage{hyperref}
\usepackage{amsthm,amsfonts,amsmath,amscd,amssymb}
\usepackage[capitalise,nameinlink]{cleveref}
\usepackage{float}
\usepackage[all]{xy}
\usepackage{mathrsfs}
\usepackage{graphics,graphicx}
\usepackage{caption}
\usepackage{subcaption}
\usepackage{comment}
\usepackage{array}
\newcolumntype{P}[1]{>{\centering\arraybackslash}p{#1}}
\newcolumntype{M}[1]{>{\centering\arraybackslash}m{#1}}
\usepackage[mathscr]{euscript}

\DeclareMathAlphabet\EuScript{U}{eus}{m}{n}
\SetMathAlphabet\EuScript{bold}{U}{eus}{b}{n}

\usepackage{epigraph}
\let\oldmarginpar\marginpar
\renewcommand\marginpar[1]{\-\oldmarginpar[\raggedleft\footnotesize #1]%
	{\raggedright\footnotesize #1}}
\theoremstyle{plain}
\newtheorem{theorem}{Theorem}[section]
\newtheorem{lemma}[theorem]{Lemma}

\newtheorem*{theorem*}{Theorem}
\newtheorem*{corollary*}{Corollary}

\newtheorem{goal}[theorem]{Goal}
\newtheorem{prop}[theorem]{Proposition}

\newtheorem{cor}[theorem]{Corollary}

\newtheorem{conj}[theorem]{Conjecture}
\theoremstyle{definition}

\newtheorem{definition}[theorem]{Definition}

\newtheorem{remark}[theorem]{Remark}

\numberwithin{equation}{section}

\crefname{theorem}{Theorem}{Theorems}
\crefname{question}{Question}{Questions}

\makeatletter
\g@addto@macro\normalsize{%
  \setlength{\abovedisplayskip}{6pt}%
  \setlength{\belowdisplayskip}{6pt}%
  \setlength{\abovedisplayshortskip}{3pt}%
  \setlength{\belowdisplayshortskip}{3pt}%
}
\makeatother

\renewcommand{\S}{\mathcal{S}}

\newcommand{\D}{\mathbb{D}}
\newcommand{\N}{\mathbb{N}}
\newcommand{\Z}{\mathbb{Z}}

\newcommand{\R}{\mathbb{R}}

\newcommand{\SA}{\mathscr{A}}

\newcommand{\SF}{\mathscr{F}}

\newcommand{\A}{\mathcal{A}}
\newcommand{\La}{\Lambda}

\newcommand{\la}{\lambda}

\newcommand{\dd}{\partial}

\newcommand{\sse}{\subseteq}
\newcommand{\lr}{\longrightarrow}

\newcommand{\GL}{\operatorname{GL}}

\newcommand{\Sh}{\operatorname{Sh}}

\newcommand{\st}{\text{st}}

\newcounter{daggerfootnote}

\newcommand{\bA}{\mathbb{A}}

\newcommand{\bG}{\mathbb{G}}

\newcommand{\M}{\mathfrak{M}}

\newcommand{\Loc}{\mathrm{Loc}}

\usepackage{dsfont}
\allowdisplaybreaks

\newcommand{\id}{\mathrm{id}}
\def \vertbar [#1](#2,#3,#4){
    \draw [#1] (#2,#3) -- (#2,#4);
    \draw [fill=white] (#2,#3) circle [radius=0.1];
    \draw [fill=black] (#2,#4) circle [radius=0.1];
}

\usepackage{mdframed}
\usepackage{lipsum}
\newtheorem*{theo}{Main Theorem}

\mdfdefinestyle{style1}{leftmargin=5cm,rightmargin=5cm,skipbelow=1cm,skipabove=1cm,
   innertopmargin=0pt}

\usepackage{varwidth}
\definecolor{DarkGreen}{RGB}{1,50,32}

\newcommand{\G}[0]{\mathbb{G}}

\newcommand{\im}{\operatorname{im}}

\newcommand{\sbs}{\hfill\rule{1.2ex}{1.2ex}}

\title{Decomposable concordances induce closed immersions}
\author{Roger Casals}
\date{}
\begin{document}

\maketitle
\begin{abstract} We show that a decomposable Lagrangian concordance between two Legendrian links must induce a closed immersion between the sheaf moduli of its ends. As the main application of this criterion, we prove a conjecture of Cornwell-Ng-Sivek, showing that their Lagrangian concordance is not decomposable. This establishes the first known example of a regular Lagrangian concordance between Legendrian knots that is not decomposable. In addition, we show that even the simpler Lagrangian concordance, satelliting via the Reeb 2-copy, is not decomposable.
\end{abstract}


\section{Introduction}

The object of this article will be to show that a decomposable Lagrangian concordance induces a closed immersion between the sheaf moduli of its Legendrian ends and illustrate that, in general, this does not occur for all regular Lagrangian concordances. Specifically, we show that the Cornwell-Ng-Sivek concordance is not decomposable, as initially conjectured 12 years ago.  This provides the first known example of a regular Lagrangian concordance that is not decomposable, which is the Lagrangian analogue of establishing the existence of a non-ribbon concordance in the context of Lagrangian surfaces in Weinstein 4-manifolds. The article also includes a second application, showing that the Reeb 2-copy variant of their Lagrangian concordance is not decomposable either.


\subsection{Scientific context} The study and understanding of embedded exact Lagrangian cobordisms between Legendrian links is a pillar of low-dimensional contact and symplectic topology. It is a rather broad landscape, with questions that range from the classification of Lagrangian fillings to the relation between such Lagrangian cobordisms and contact and symplectic invariants. See e.g.~\cite{BourgeoisSabloffTraynor15,Chantraine10,CRGG20,Chekanov,EHK,Etnyre05,etnyre2024lagrangianrealizationsribboncobordisms,EtnyreNg22_Survey,guadagni2022legendriansatellitesdecomposablecobordisms,HenryRutherford15,li2024lagrangiancobordismfunctormicrolocal,YuPan2,rizell2026nonregularlagrangianconcordanceslagrangian,Traynor01}. In the context of Lagrangian cobordisms, the class of decomposable Lagrangian cobordisms, cf.~\cite[Section 6]{EHK} or \cite[Section 3]{CornwellNgSivek16}, has a role akin to that of ribbon cobordisms in low-dimensional smooth topology.

\noindent In 2014, the article \cite{CornwellNgSivek16} brought forth the first candidate for a non-decomposable Lagrangian concordance between Legendrian knots. The construction uses a Lagrangian ribbon concordance from the max-tb Legendrian unknot to the max-tb Legendrian $m(9_{46})$, and then satellites this concordance with a specific Legendrian pattern based on an appropriately framed Legendrian Whitehead double. It was conjectured in \cite[Conjecture 3.3]{CornwellNgSivek16} that the resulting Lagrangian concordance $L$ is not decomposable. Their construction is particularly appealing, as there are no immediate reasons why the resulting Lagrangian concordance should not be decomposable:
\begin{enumerate}
    \item By \cite[Section 3]{KAWAUCHI1980321}, see also \cite{Baker_2015}, smooth satellite operations preserve ribbon concordance. Thus $L$ is smoothly ribbon, i.e.~the smooth ribbon condition does not obstruct $L$ from being Lagrangian decomposable.\\

    \item By \cite[Theorem 1.11]{breen2024regularlysliceimpliesoncestably}, the Lagrangian concordance $L$ is regular, i.e.~compatible with the natural ambient Weinstein structure. Decomposable Lagrangian cobordisms are necessarily regular. Therefore, there are no regularity obstructions for $L$ to be Lagrangian decomposable.
\end{enumerate}
\noindent Note that it is possible to construct interesting non-regular Lagrangian concordances, cf.~\cite{rizell2026nonregularlagrangianconcordanceslagrangian}. That said, the focus is to address the subtler phenomenon presented in \cite[Conjecture 3.3]{CornwellNgSivek16}, and it does not involve stabilizations of any Legendrian links.

A relevant aspect of \cite[Conjecture 3.3]{CornwellNgSivek16}, i.e.~of deciding whether $L$ is Lagrangian decomposable or not, is that establishing it either way requires a new idea. If one is to prove that $L$ is Lagrangian decomposable, then one must develop techniques to manipulate Lagrangian cobordisms so as to simplify them to decomposable ones. If one is to obstruct $L$ from being Lagrangian decomposable, one must find a symplectic geometric property satisfied by decomposable Lagrangian concordances that is not satisfied by a general Lagrangian concordance. The first important conceptual contribution of this article is that such a property exists and can be used to detect non-decomposable Lagrangian concordances, cf.~\cref{lem:decomposable_closed}.


\subsection{Main results} Let $L:\La_-\prec\La_+$ be an embedded exact Lagrangian concordance. By \cite[Theorem 1.1]{li2024lagrangiancobordismfunctormicrolocal}, such a Lagrangian concordance induces a morphism
$$\Phi_L:\M_r(\La_-)\lr\M_r(\La_+)$$
between the derived stacks of microlocal rank $r$ pseudoperfect sheaves on $\R^2$ with singular support on $\La_\pm$, cf.~\Cref{ssec:initial_ingredients}. It is known to be an open immersion, see e.g.~\cite[Theorem 1.1]{li2024lagrangiancobordismfunctormicrolocal} or \cite[Prop.~B.26]{CasalsLi}. A central contribution in this article is the observation that this induced morphism is also a closed immersion if $L$ is decomposable. Specifically, this criterion reads as follows:

\begin{theorem}
\label{lem:decomposable_closed} Let $L\sse(\R^4,\la_\st)$ be a decomposable Lagrangian concordance $L:\La_-\prec \La_+$ and $r\in\N$. Then the induced morphism
\begin{equation*}
\Phi_L:\M_r(\La_-)\lr \M_r(\La_+)
\end{equation*}
is a closed immersion.\sbs
\end{theorem}

\noindent \Cref{lem:decomposable_closed} is summarized in the title of this article. While it is conceptually useful to know that \Cref{lem:decomposable_closed} holds, it only becomes particularly valuable once it can be effectively used to detect that non-decomposable Lagrangian concordances exist. The second part of this article focuses on such a task. Specifically, we focus on establishing \cite[Conjecture 3.3]{CornwellNgSivek16}, as we now describe in detail.

\begin{center}
	\begin{figure}[H]
		\centering
		\includegraphics[scale=0.8]{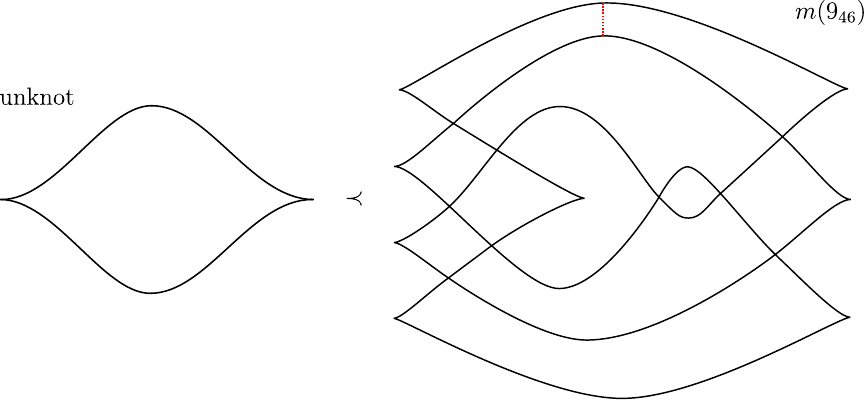}
		\caption{(Left) A front for the max-tb Legendrian unknot. (Right) A front for the max-tb Legendrian $\La_{m(9_{46})}$. A decomposable Lagrangian concordance can be obtained from the left to the right by pinching on the vertical red line and filling a resulting component with a minimum cobordism.}
		\label{fig:Front_SatelliteConcordance1}
	\end{figure}
\end{center}

Let us consider the max-tb Legendrian representative $\La_{m(9_{46})}$ of the smooth knot $m(9_{46})$. A front for it is depicted in \Cref{fig:Front_SatelliteConcordance1} (right). Such a Legendrian knot bounds a decomposable Lagrangian disk. Specifically, by performing a pinch move on the red vertical arc depicted in \Cref{fig:Front_SatelliteConcordance1} (right) we obtain an orientable embedded exact Lagrangian cobordism down to the 2-component standard Legendrian unlink. One of these components can be filled with a Lagrangian 2-disk and we obtain a decomposable Lagrangian concordance $\ell:\La_{\st}\prec\La_{m(9_{46})}$, from the max-tb unknot $\La_\st$ to $\La_{m(9_{46})}$.

\begin{center}
	\begin{figure}[H]
		\centering
		\includegraphics[scale=0.8]{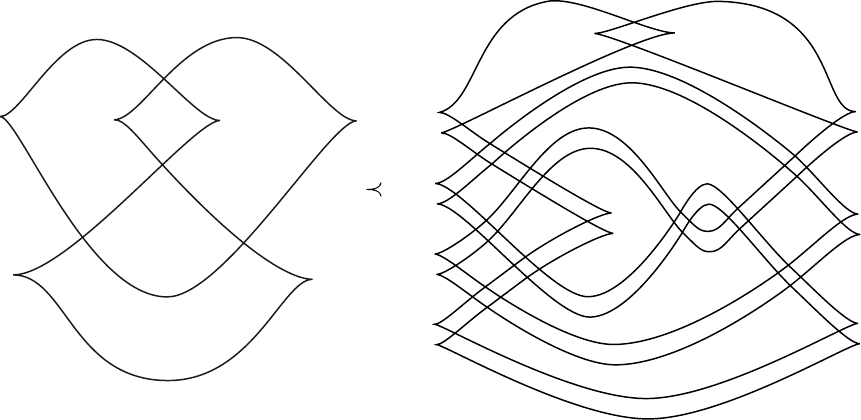}
		\caption{(Left) A front for the max-tb Legendrian right-handed trefoil. (Right) A front for the Legendrian Whitehead double $W(\La_{m(9_{46})})$. It is proven in \cite{CornwellNgSivek16} that there is an embedded exact Lagrangian concordance from the former to the latter. \Cref{thm:main} proves that it is not Lagrangian decomposable.}
		\label{fig:Front_SatelliteConcordance2}
	\end{figure}
\end{center}

\noindent The main construction relevant to \cite[Conjecture 3.3]{CornwellNgSivek16} satellites the Lagrangian concordance $\ell:\La_{\st}\prec\La_{m(9_{46})}$ via a Legendrian Whitehead double pattern, cf.~\cite[Section 3]{CornwellNgSivek16}. It is proven in \cite{CornwellNgSivek16} that the result of such a satellite is an embedded exact Lagrangian concordance $L:W(\La_\st)\prec W(\La_{m(9_{46})})$ from the Legendrian Whitehead double $W(\La_\st)$ of $\La_\st$ to the Legendrian Whitehead double $W(\La_{m(9_{46})})$ of $\La_{m(9_{46})}$, cf.~\cite[Theorem 2.4]{CornwellNgSivek16}. It is explained in \cite[Figure 8]{CornwellNgSivek16} that $W(\La_\st)$ is Legendrian isotopic to the unique max-tb representative $\La_{3_1}$ of the right-handed trefoil knot. Thus we denote the concordance by $L:\La_{3_1}\prec W(\La_{m(9_{46})})$. \Cref{fig:Front_SatelliteConcordance2} depicts Legendrian fronts for the concave end $\La_{3_1}$ of this Lagrangian concordance, on its left, and for its convex end $W(\La_{m(9_{46})})$, on its right. We refer to this Lagrangian concordance as the Cornwell-Ng-Sivek concordance. As said above, it is conjectured in \cite[Conjecture 3.3]{CornwellNgSivek16} that this Lagrangian concordance is not decomposable. The second contribution of this article is to establish such a result:

\begin{cor}\label{thm:main} The Cornwell-Ng-Sivek Lagrangian concordance
$$L:\La_{3_1}\prec W(\La_{m(9_{46})})$$
is a regular Lagrangian concordance which is not decomposable.\sbs
\end{cor}

\Cref{thm:main} is proven by using \Cref{lem:decomposable_closed} and developing a better understanding of the specific morphism $\M_1(\La_-)\lr\M_1(\La_+)$, and its domain and codomain, in the case of the Lagrangian concordance at hand. We emphasize that \Cref{thm:main} also settles the longstanding question of whether regularity of a Lagrangian surface is equivalent to decomposability. The conclusion is that a decomposable Lagrangian cobordism must be regular but there are regular Lagrangian cobordisms that are not decomposable.

\begin{remark} The max-tb Legendrian $m(9_{46})$ admits at least two embedded Lagrangian 2-disk fillings, see e.g.~\cite[Example 2.8]{CornwellNgSivek16}. Specifically, one can pinch at the red vertical line in \Cref{fig:Front_SatelliteConcordance1} (right) or at the corresponding vertical line between the corresponding strands at the bottom of the front. It is not specified which of these two Lagrangian 2-disks is being used for \cite[Conjecture 3.3]{CornwellNgSivek16}. Our proof of \Cref{thm:main} shows that either choice leads to a non-decomposable concordance, i.e.~\Cref{thm:main} applies to both such Lagrangian concordances.\sbs
\end{remark}

Finally, \Cref{lem:decomposable_closed} can be used to detect many other non-decomposable Lagrangian concordances. Let us illustrate this by considering a construction even simpler than \cite[Conjecture 3.3]{CornwellNgSivek16}, namely, by taking the Reeb 2-copy of one of the two Lagrangian concordances $\La_\st\prec\La_{m(9_{46})}$. The Reeb 2-copy of the Legendrian max-tb unknot is the standard Legendrian Hopf link $\La_H\simeq\La^{(2)}_\st$, and we denote the Reeb 2-copy of $\La_{m(9_{46})}$ by $\La_{m(9_{46})}^{(2)}$. A similar argument to that establishing \Cref{thm:main} also shows the following result:

\begin{cor}\label{thm:main2} The Reeb 2-copy
$$L^{(2)}:\La_H\prec \La_{m(9_{46})}^{(2)}$$
of either of the two specified Lagrangian concordances $\ell_\pm:\La_\st\prec\La_{m(9_{46})}$, obtained by puncturing either of the two specified Lagrangian 2-disks of $\La_{m(9_{46})}$, is not decomposable.\sbs
\end{cor}

\noindent \Cref{thm:main2} helps illustrate that the clasp in the Whitehead double is not the sole cause, in a sense, of the non-decomposability of the Cornwell-Ng-Sivek concordance. Rather, non-decomposability already appears when considering the Reeb 2-copy, even before introducing the clasp. It is also valuable to have both \Cref{thm:main} and \Cref{thm:main2}, as \Cref{thm:main} shows the existence of non-decomposable regular Lagrangian concordances between Legendrian knots, not just Legendrian links.\\


\noindent {\bf Acknowledgements}. I am thankful to J.~Etnyre, L.~Ng, B.~Tosun and L.~Traynor for organizing the NSF CBMS Conference ``Legendrian Links and the Microlocal Theory of Sheaves''. The results of this article stem from working on the notes for that conference. The author is supported by the National Science Foundation under the grant DMS-2505760.\sbs

\vspace{0.2cm}
\noindent {\bf Conventions}. Legendrian links are always considered in the standard contact $(\R^3,\xi_\st)$. Lagrangian cobordisms between Legendrian links are always considered in the symplectization $(\R^4,\la_\st)$ of $(\R^3,\xi_\st)$. All Legendrian links appearing in this note have rotation zero and, thus, the sheaf-theoretic invariants can have values in $\Z$-graded cochain complexes. In this article the word ruling always means a $\Z$-graded normal ruling.\sbs


\section{From decomposable concordances to closed immersions}

The goal of this section is to establish \Cref{lem:decomposable_closed}, showing that the sheaf quantization of a decomposable Lagrangian concordance must necessarily induce a closed immersion between the sheaf moduli of its ends. In conjunction with the results from \Cref{sec:nondecomposability_CNS}, this will establish \Cref{thm:main} and \Cref{thm:main2}.

\subsection{The initial ingredients}\label{ssec:initial_ingredients}
Let $\La\sse(\R^3,\xi_\st)$ be a Legendrian link and $\pi(\La)\sse\R^2$ a Legendrian front for $\La$. Consider the dg derived category $\Sh_\La(\R^2)$ of sheaves on $\R^2$ with singular support along $\La$ and coefficients in the dg category of cochain complexes of $k$-vector spaces, cf.~\cite[Chapter V]{KashiwaraSchapira_Book}, \cite[Section 3.1]{STZ2017} or \cite[Section 2.8]{CW}. Let us consider the smooth dg derived subcategory $\Sh^{c}_\La(\R^2)_0$ of categorically compact objects whose stalk in the unbounded region of $\R^2\setminus \pi(\La)$ is acyclic, cf.~\cite[Section 3]{Nadler16}.

We denote by $\M(\La)$ the derived moduli stack of pseudoperfect modules over $\Sh^{c}_\La(\R^2)_0$, as in \cite[Section 3]{ToenVaquie07}, and by $\M_r(\La)$ the derived stack that parametrizes pseudoperfect sheaves of a fixed microlocal rank $r\in\N$, see e.g.~\cite[Section 2]{casals2024positivemicrolocalholonomiesglobally}. In particular, it follows from \cite[Theorem 3.7]{GKS_Quantization}, see also \cite[Theorem 4.9]{STZ2017}, that $\M_r(\La)$ is a Legendrian invariant.

\subsection{Decomposable Lagrangian concordances induce closed immersions}\label{ssec:property_decomposable}

By \cite[Theorem 1.1]{LiCobordism1}, cf.~also \cite[Theorem 1.1]{li2024lagrangiancobordismfunctormicrolocal}, a Lagrangian cobordism $L:\La_-\prec\La_+$, with the appropriate grading and relative spin hypothesis, induces an algebraic morphism

\begin{equation}\label{eq:map_cobordism}
\Phi_L:\M_r(\La_-)\times_{\Loc_r(\La_-)}\Loc_r(L)\lr \M_r(\La_+)
\end{equation}

\noindent of derived stacks, where $\Loc$ are the corresponding moduli stacks for the dg categories of perfect local systems on $\La_-$ and $L$. In general, $\Phi_L$ is an open immersion of derived stacks. Indeed, \cite[Theorem 1.1]{LiCobordism1} shows that the functor induced on the proper dg categories of sheaves is fully faithful, which implies that the map between the corresponding stacks is an \'etale monomorphism. Since these stacks are of finite presentation, an \'etale monomorphism must be an open immersion. See also~\cite[Prop.~B.26]{CasalsLi} and references therein, which already proved such a statement for Lagrangian cobordisms with $\La_-=\emptyset$.

The morphism $\Phi_L$ is not, in general, a closed immersion. To wit, for any genus $g$ embedded exact Lagrangian filling of a Legendrian knot $\La_+$ and $r=1$, the morphism $\Phi_L$ is an open immersion
$$\Phi_L:\Loc_1(L)\simeq \G_m^{2g}\times B\G_m\hookrightarrow\M_1(\La_+)$$
for which the image of the domain $B\G_m\times\Loc_1(L)\cong B\G_m\times \G_m^{2g(L)}$ is not Zariski closed in $\M_1(\La_+)$ for almost any $\La_+$. Specifically, in these cases, the image is dense but it does not coincide with the entire moduli $\M_1(\La_+)$. This already occurs if we consider $\La_+$ to be the max-tb Legendrian trefoil knot and $L$ any of its embedded exact Lagrangian fillings from \cite{EHK} or \cite{YuPan}.

\noindent Let us now show that if $L$ is a decomposable Lagrangian concordance then $\Phi_L$ is a closed immersion, i.e.~prove \Cref{lem:decomposable_closed}, restated here for convenience:

\begin{theorem}
\label{lem:decomposable_closed2} Let $L\sse(\R^4,\la_\st)$ be a decomposable Lagrangian concordance $L:\La_-\prec \La_+$ and $r\in\N$. Then the induced morphism $\Phi_L:\M_r(\La_-)\lr \M_r(\La_+)$ is a closed immersion.
\end{theorem}

\begin{proof} By definition, cf.~\cite[Section 6]{EHK}, a decomposable Lagrangian cobordism, scanned bottom to top, consists of the following three types of embedded exact Lagrangian cobordisms:

\begin{enumerate}
    \item The trace of a Legendrian isotopy, cf.~\cite[Section 6.1]{EHK}. In this piece of the Lagrangian cobordism, the concave and convex Legendrian ends are Legendrian isotopic.\\
    
    \item A Lagrangian index-0 handle attachment, i.e.~ a minimum cobordism, cf.~\cite[Section 6.4]{EHK} or \cite[Section 2.1]{CasalsNg}. In such a piece of the Lagrangian cobordism, the convex Legendrian end is Legendrian isotopic to the Legendrian split union of the concave end with a max-tb Legendrian unknot.\\
    
    \item A Lagrangian index-1 handle attachment, i.e.~ a saddle cobordism, cf.~\cite[Section 6.5]{EHK} or \cite[Section 2.1]{CasalsNg}. A front for the concave Legendrian end of such a Lagrangian cobordism can be obtained by performing a pinch move to a front for the convex Legendrian end, cf.~\cite[Figure 11]{EHK}.\\
\end{enumerate}

\noindent Let us argue that for each of these three types of Lagrangian cobordisms, the corresponding induced map is a closed immersion under the given hypotheses. To ease notation, we denote the given decomposable Lagrangian concordance $L^{tot}:\La^{tot}_-\prec \La^{tot}_+$ and reuse the notation $L:\La_-\prec\La_+$ for each type of elementary Lagrangian cobordism.

First, if $L:\La_-\prec\La_+$ is the trace of a Legendrian isotopy, then $L$ is an invertible Lagrangian concordance and thus $\Phi_L$ is an isomorphism. In particular, it is a closed immersion.

Second, if $L:\La_-\prec\La_+$ is a minimum cobordism, then $L$ consists of a disjoint union of the Lagrangian cylinder over $\La_-$ and a Lagrangian 2-disk $\D^2$. In particular, the convex end $\La_+=\La_-\sqcup \La_\st$ is a disjoint unknotted union of $\La_-$ with the max-tb standard Legendrian unknot $\La_\st$. Since the Lagrangian $\D^2$-component of $L$ is unknotted from the cylinder over $\La_-$, the domain of $\Phi_L$ simplifies to
\begin{equation}\label{eq:minimum_cob1}
\M_r(\La_-)\times_{\Loc_r(\La_-)}\Loc_r(L)\cong\M_r(\La_-)\times \Loc_r(\D^2)
\end{equation}
To conclude, note that $\Loc_r(\D^2)\cong B\GL_{r}$ and $\M_r(\La_\st)\cong B\GL_{r}$ and the map induced by such Lagrangian 2-disk $\D^2$ as a Lagrangian filling of $\La_\st$ precisely induces an isomorphism between them. Thus, in combination with \eqref{eq:minimum_cob1}, the induced map
$$\Phi_L:\M_r(\La_-)\times_{\Loc_r(\La_-)}B\GL_{r}\lr \M_r(\La_-\sqcup\La_\st)$$
is equivalent to
$$\Phi_L:\M_r(\La_-)\times \Loc_r(\D^2)\lr\M_r(\La_-)\times B\GL_{r}$$
where we have used the isomorphism $\M_r(\La_-\sqcup\La_\st)\cong\M_r(\La_-)\times B\GL_{r}$, as $\La_\st$ is Legendrian unknotted from $\La_-$, and $\Phi_L$ is an isomorphism when restricted to either factor. In conclusion, $\Phi_L$ is also an isomorphism of derived stacks for the minimum cobordism $L$, and so it is a closed immersion.

\begin{center}
	\begin{figure}[H]
		\centering
		\includegraphics[scale=0.65]{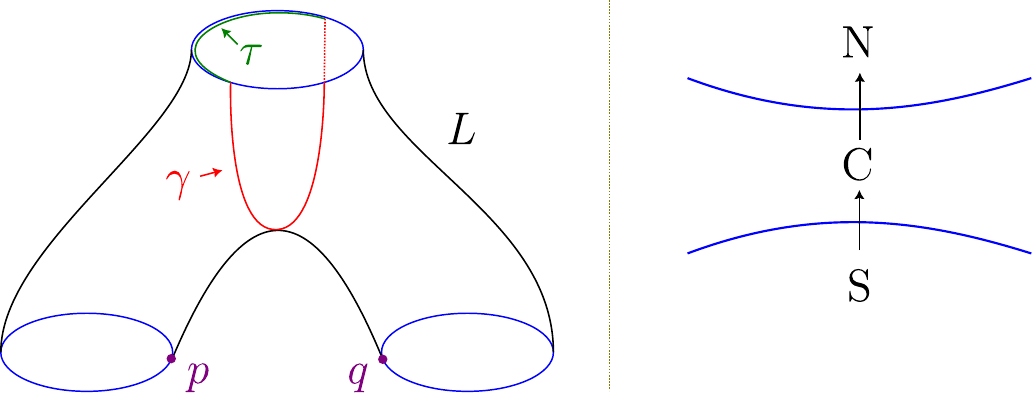}
		\caption{(Left) An elementary saddle cobordism $L:\La_-\prec\La_+$. (Right) The local region for the front of $\La_+$ where a pinching move along its unique Reeb chord results in the elementary saddle cobordism, cf.~\cite[Section 6.5]{EHK}. Here $S,C$ and $N$ are the stalks of a sheaf in $\R^2$ with singular support on $\La_+$, and the arrows indicate corresponding maps between the stalks.}
		\label{fig:SaddleCob}
	\end{figure}
\end{center}

Third, if $L:\La_-\prec\La_+$ is a saddle cobordism, then the relevant component of $L$ is smoothly isomorphic to a pair of pants, the other components being cylinders. The number of components of $\La_+$ equals the number of components of $\La_-$ plus or minus one. For a Lagrangian saddle cobordism, the map $\Phi_L$ is studied in detail in \cite[Section 4.1.1]{LiCobordism1}. In particular, we note that $\Phi_L$ is {\it not} in general a closed immersion. Its image $\mbox{im}(\Phi_L)\sse\M(\La_+)$ can be described as follows. Consider the pinch move from $\La_+$ to $\La_-$ that induces the saddle cobordism $L$: at the local model for $\La_+$ near where the pinch move is to occur, we have two strands of $\La_+$, as depicted in \Cref{fig:SaddleCob} (right). Let us refer to the stalks of a sheaf in $\M(\La_+)$ as $S,C,N$, where $S$ is the stalk for the region south of the bottom strand, $C$ is the stalk in the central region between the two strands and $N$ is the region north of the top strand. The singular support condition gives two morphisms
$$\beta:S\lr C,\quad \alpha:C\lr N.$$

\noindent Then $\mbox{im}(\Phi_L)\sse\M(\La_+)$ is characterized as the substack of sheaves in $\M(\La_+)$ such that $\alpha\circ\beta:S\lr N$ is an equivalence. That this is necessary follows from the fact that, after the pinch move, the regions with stalks $S$ and $N$ become the same region, see also \cite[Section 4.1.1]{LiCobordism1}. Without any further hypotheses, $\mbox{im}(\Phi_L)\sse\M(\La_+)$ is not always Zariski closed.

\noindent Let us now use the hypothesis that $L^{tot}$ is a decomposable Lagrangian concordance, which we do for the first and only time. Specifically, we only need to use that $L^{tot}$ is a smooth ribbon concordance. Since $L^{tot}$ is a smooth ribbon concordance, any of its elementary saddle cobordisms must be such that at least one of the components to which the 1-handle is attached bounds a disk in $L^{tot}$. That is, if the core of the Lagrangian 1-handle of the saddle cobordism $L\sse L^{tot}$ intersects $\La_-$ at points $p,q\in\La_-$, cf.~\Cref{fig:SaddleCob} (left), then either the component of $\La_-$ that contains $p$ bounds a Lagrangian 2-disk contained in the part of $L^{tot}$ below and including $L$, or the component of $\La_-$ that contains $q$ does. Therefore, thanks to the hypothesis that $L$ is a decomposable Lagrangian concordance, we can and do assume that all the elementary saddle cobordisms that occur in $L^{tot}$ are such that the ends of the Lagrangian 1-handles are attached to different components and at least one of them bounds a Lagrangian disk in $L^{tot}$.\\

\noindent Thus, we need to study the case that $L:\La_-\prec\La_+$ is a Lagrangian saddle cobordism $L\sse L^{tot}$ such that the ends $p,q\in\La_-$ of the core of the Lagrangian 1-handle are attached to different connected components $\La_p,\La_q$ of $\La_-$ and $\La_p$, say, bounds a Lagrangian disk in $L^{tot}$, cf.~\Cref{fig:SaddleCob} (left). Consider the cocore $\gamma\sse L$ of the saddle cobordism $L$ and the closure $\tau\sse\La_+$ of the arc in the component relevant to $\La_+\setminus\dd\gamma$ for which the concatenation $\gamma\circ\tau$ is homotopic to $\La_p$, as seen within the elementary cobordism $L$. The cocore $\gamma$ and $\tau$ are depicted in red and green respectively in \Cref{fig:SaddleCob} (left). Denote by $\dd_\pm\gamma\sse\La_+$ the endpoints of the cocore, with $\dd_-\gamma$ and $\dd_+\gamma$ respectively lying in the bottom and top strands depicted in \Cref{fig:SaddleCob} (right).\\

Recall that above we described a sheaf in $\M(\La_+)$ near the pinching region by two morphisms
$$\beta:S\lr C,\quad \alpha:C\lr N$$

\noindent and $\mbox{im}(\Phi_L)\sse\M(\La_+)$ consists of sheaves in $\M(\La_+)$ such that $\alpha\circ\beta:S\lr N$ is an equivalence. In these terms, the microlocal stalk $\mu_-(\SF)$ at $\dd_-\gamma\in\La_+$ of a sheaf $\SF\in\M(\La_+)$ is $\mbox{cofib}(\beta)$, while the microlocal stalk $\mu_+(\SF)$ at $\dd_+\gamma\in\La_+$ is $\mbox{fib}(\alpha)$. The two useful observations now are:

\begin{enumerate}
    \item Consider the natural map
    \begin{equation}
        \vartheta:\mbox{fib}(\alpha)\lr\mbox{cofib}(\beta)
    \end{equation}
    factoring through $C$. Since, by construction, $\mbox{cofib}(\vartheta)\simeq \mbox{cofib}(\alpha\circ\beta)$, $\mbox{im}(\Phi_L)$ can be equivalently described as those sheaves for which $\vartheta$ is an equivalence.\\

    \item We claim that the additional fact that $\La_p$ bounds a Lagrangian disk in $L^{tot}$ forces the equivalence $\vartheta$ to be equal to the inverse of the microlocal parallel transport
    $$T_\gamma:\mu_+(\SF)\lr\mu_-(\SF)$$
    along $\gamma$.\footnote{Here we are implicitly using that $L$ admits a sheaf quantization, which it does, cf.~\cite{LiCobordism1}.} This is the key passage to a closed condition. Indeed, the claim holds because the concatenation $\gamma\circ\tau$ is homotopic to $\La_p$ in $L$. Therefore, since $\La_p$ bounds a Lagrangian 2-disk in $L^{tot}$, the microlocal parallel transport along $\gamma\circ\tau$ must be trivial. In consequence $T_\gamma^{-1}\circ\vartheta=\id_{\mu_+(\SF)}$, which is a closed condition.
\end{enumerate}

\noindent Therefore, under the additional assumption that a saddle cobordism $L:\La_-\prec\La_+$ is part of a decomposable Lagrangian concordance $L^{tot}$, we conclude that $\Phi_L$ is a closed immersion. This establishes all the cases and thus the required statement follows.\end{proof}

\section{Non-decomposability of the Cornwell-Ng-Sivek concordance}\label{sec:nondecomposability_CNS}

The goal of this section is to use \Cref{lem:decomposable_closed} to show that the Cornwell-Ng-Sivek concordance is not decomposable, establishing \Cref{thm:main}. The main techniques that we use are from the microlocal theory of sheaves and the associated ruling decompositions, see e.g.~\cite{KashiwaraSchapira_Book,GKS_Quantization,Guillermou23_SheafSummary,HenryRutherford15,STZ2017}. We also establish \Cref{thm:main2} using \Cref{lem:decomposable_closed}.


\subsection{A sheaf moduli for the concave end}

Let $\La_{3_1}$ denote the unique max-tb representative of the Legendrian right-handed trefoil. Since $\La_{3_1}$ is Legendrian isotopic to the rainbow closure of the 2-stranded braid word $\sigma_1^3$, the derived stack $\M_1(\La_{3_1})$ can be computed directly via the methods of \cite[Section 3.4]{STZ2017}, see also \cite[Section 4.1]{CW}. Alternatively, it can be deduced by using the Legendrian dg algebra of $\La_{3_1}$ in combination with \cite{NRSSZ}. That said, since $\La_{3_1}$ is naturally presented as the Legendrian Whitehead double of the max-tb unknot $\La_\st$ in \cite[Conjecture 3.3]{CornwellNgSivek16}, and thus its front is given by \Cref{fig:Front_SatelliteConcordance2} (left), it might be enlightening to record the computation of $\M_1(\La_{3_1})$ directly in this front, as follows.

\begin{lemma}\label{lem:sheaves_trefoil} Let $\La_{3_1}$ denote the unique max-tb representative of the Legendrian right-handed trefoil. Then there is an isomorphism
\begin{equation}\label{eq:moduli_trefoil}
\M_1(\La_{3_1}) \cong \{(x,y,z)\in\mathbb{A}^3: xyz+x+z=1\}\times B\bG_m.
\end{equation}
\end{lemma}

\begin{proof} Consider the front in \Cref{fig:Front_SatelliteConcordance2} (left) and the unique integral Maslov potential which has potential $0$ at the lower stand of the front. Throughout the argument we will now use the local models discussed in \cite[Section 3.4]{STZ2017}, or see also \cite[Section 4]{ChantraineNgSivek} or \cite[Section 4.1]{CW}. That is, we present the sheaves on $\R^2$ with singular support on $\La_{3_1}$ as certain constructible sheaves on $\R^2$ and, due to the exodromy equivalence \cite[Theorem A.9.3]{HA}, we use that those constructible sheaves themselves can be described through the representation theory of the poset associated to the stratification of $\R^2$ given by the front of $\La_{3_1}$.\\

\begin{center}
	\begin{figure}[H]
		\centering
		\includegraphics[scale=1]{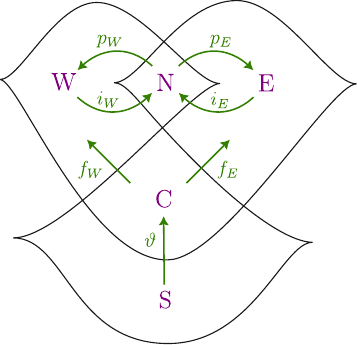}
		\caption{A front for the max-tb Legendrian right-handed trefoil $\La_{3_1}$, as the Legendrian Whitehead double of the max-tb Legendrian unknot. We have highlighted the stalks of a sheaf $\SF$ singularly supported on the cone of $\La_{3_1}$, as well as the necessary maps to describe such a sheaf $\SF$.}
		\label{fig:Front_31_sheaves}
	\end{figure}
\end{center}

In this particular case, the exodromy equivalence along with the singular support conditions implies that $\M_1(\La_{3_1})$ is the moduli stack of dg representations of the diagram drawn in \Cref{fig:Front_31_sheaves} with certain additional constraints, according to \cite[Section 3.4]{STZ2017}, see also \cite[Section 4]{ChantraineNgSivek}. In particular, the square diagram of arrows at a crossing must be a homotopy Cartesian square. Recall that this moduli stack should also take into account the natural actions of the corresponding change-of-basis gauge group at each vertex of the quiver. The Maslov potential immediately imposes the following constraints:

\begin{enumerate}
    \item $S$ is a 1-dimensional $k$-vector space concentrated in degree 0.
    \item $E,W$ are both 1-dimensional $k$-vector spaces concentrated in cohomological degree -1.
    \item $N$ is a 2-dimensional $k$-vector space concentrated in cohomological degree -1.
\end{enumerate}

\noindent The stalk $C$ depends on the sheaf: it is a 2-term complex but whether it is acyclic or not depends on the point of $\M_1(\La_{3_1})$. The singular support conditions also require
\begin{equation}\label{eq:trefoil_sheaves2}
N\cong\ker(p_W)\oplus\ker(p_E)
\end{equation}
because of the crossing in the top part of $N$, as well as
\begin{equation}\label{eq:trefoil_sheaves1}
p_W\circ i_W=\mbox{id}_W,\qquad p_E\circ i_E=\mbox{id}_E
\end{equation}
\noindent because of the right and left cusps at $N$. By the constraint \eqref{eq:trefoil_sheaves2}, we can use the $\GL_2(k)$ local gauge at $N$ in order to ensure that the projections $p_E$ and $p_W$ are written as
\begin{equation}\label{eq:trefoil_sheaves3}
p_W=(1,0),\qquad p_E=(0,1),
\end{equation}
and thus \eqref{eq:trefoil_sheaves1} implies that the inclusions $i_W,i_E$ can be written as
\begin{equation}\label{eq:trefoil_sheaves4}
i_W=\binom{1}{-x},\qquad i_E=\binom{y}{1}
\end{equation}
for some variables $x,y$ with values in the base field.

Let us now impose the singular support condition at the crossing at the bottom of $N$, i.e.~at the top of $C$. Specifically, the following homotopy-commutative square must be a homotopy pullback:
\begin{center}
\begin{tikzcd}
    & N & \\
    W \arrow[ur, "i_W"] & & E \arrow[ul, "i_E"'] \\
    & \arrow[ul, "f_W"] C \arrow[ur, "f_E"'] &
\end{tikzcd}
\end{center}
Therefore, $C$ must be homotopy equivalent to the 2-term complex
\begin{equation}\label{eq:trefoil_sheaves5}
C\simeq \mbox{fib}(W\oplus E\stackrel{(i_W,-i_E)}{\lr}N)\simeq \mbox{Cone}(W\oplus E\stackrel{(i_W,-i_E)}{\lr}N)[-1]
\end{equation}
\noindent where $W\oplus E$ is in degree -1 and $N$ in degree 0, as we use the cohomological convention. Here the maps $f_W,f_E$ become the projections from the first and second direct summands of $C_{-1}$. Note that, by \eqref{eq:trefoil_sheaves4}, the stalk \eqref{eq:trefoil_sheaves5} is acyclic if and only if $1+xy\neq0$.

Let us now impose the singular support conditions at the two remaining crossings, located left and right of $C$, equivalently at the top left and top right of $S$, see \Cref{fig:Front_31_sheaves}. The left and right crossings of $C$ respectively contain the information of null-homotopies
$$f_W\circ \vartheta\simeq 0,\qquad f_E\circ \vartheta\simeq 0,$$
which can be described by a morphism $h:S\lr W[-1]\oplus E[-1]$. By \eqref{eq:trefoil_sheaves5}, the morphism $\vartheta:S\lr C$ is given by specifying a linear map $\vartheta:S\lr N[-1]$. Note that we have the freedom to choose $\vartheta$ and $h$ up to a chain homotopy acting on the cochain complex $C$. In this case this is given by $\eta:S\lr W[-1]\oplus E[-1]$ with such $\eta$ acting via
$$\vartheta\longmapsto\vartheta+\dd_C\circ\eta,\quad h\longmapsto h+\eta.$$
\noindent Thus we can set $\eta:=-h$ and that uniquely simplifies the data of $\vartheta$ and $h$ to a unique linear map $\theta:=\vartheta-\dd_C\circ h$ from $S$ to $N[-1]$. Let us write $
\theta:S\lr N[-1]$ in coordinates as
\begin{equation}\label{eq:trefoil_sheaves52}
\theta=\binom{u}{z}
\end{equation}
for some undetermined variables $u,z$ with values in the base field. The left crossing of $C$ forces the image of $\theta$ to be linearly independent of the image of $i_W$, and the right crossing of $C$ forces the image of $\theta$ to be linearly independent with the image of $i_E$. In coordinates, using \eqref{eq:trefoil_sheaves4} and \eqref{eq:trefoil_sheaves5}, these are the conditions that
\begin{equation}\label{eq:trefoil_sheaves6}
\det\begin{pmatrix}
u & 1 \\
  z & -x
\end{pmatrix}\in k^\times,\qquad \det\begin{pmatrix}
  u & y \\
  z & 1
\end{pmatrix}\in k^\times.
\end{equation}
\noindent By using the local $\bG_m$ gauge symmetries at $W$ and $E$ respectively, we can further simplify the determinants in \eqref{eq:trefoil_sheaves6} to be equal to $\pm1$, as we choose. Therefore, after using such $\bG_m$ symmetries at $W$ and $E$, we obtain the constraints
\begin{equation}\label{eq:trefoil_sheaves7}
\det\begin{pmatrix}
  u & 1 \\
  z & -x
\end{pmatrix}=-1,\qquad \det\begin{pmatrix}
  u & y \\
  z & 1
\end{pmatrix}=1.
\end{equation}
The right equation in \eqref{eq:trefoil_sheaves7} implies $u=1+yz$, and the left equation therefore reads
\begin{equation}\label{eq:trefoil_sheaves8}
\det\begin{pmatrix}
1+yz & 1\\
z & -x
\end{pmatrix}=-1\Longleftrightarrow -x(1+yz)-z=-1\Longleftrightarrow xyz+x+z=1.
\end{equation}
Finally, after fixing $p_W=(1,0)$ and $p_E=(0,1)$, we still have the freedom of a change of basis at $W,E,S$ and a diagonal rescaling at $N$. Therefore, there is only the $\bG_m$ local symmetry at $S$ as a freedom for these local symmetries. Since we are parametrizing a 1-dimensional $k$-vector space, up to such symmetry, this freedom contributes a $B\bG_m$ factor. In conclusion, $\M_1(\La_{3_1})$ is cut out by \eqref{eq:trefoil_sheaves8} with the additional factor of $B\bG_m$, hence the statement follows.
\end{proof}

\begin{center}
	\begin{figure}[H]
		\centering
		\includegraphics[width=\textwidth]{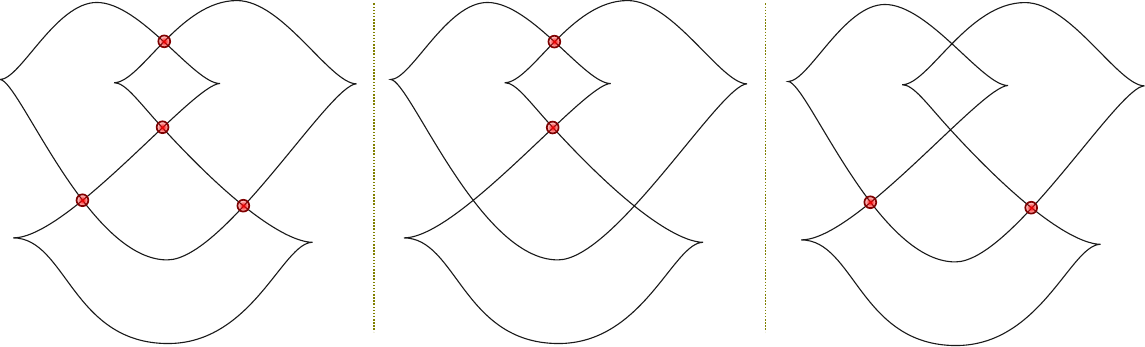}
		\caption{The three $\Z$-graded normal rulings for the max-tb Legendrian right-handed trefoil $\La_{3_1}$ presented as the Legendrian Whitehead double of the max-tb Legendrian unknot. The switches are marked with a red circle.}
		\label{fig:Front_31_ruling}
	\end{figure}
\end{center}

Two comments on \Cref{lem:sheaves_trefoil} are in order:\\

\noindent (1) The ruling decomposition for \eqref{eq:moduli_trefoil}, cf.~\cite{HenryRutherford15}, is of the form
\begin{equation}\label{eq:ruling_stratification}
\{(x,y,z)\in\mathbb{A}^3: xyz+x+z=1\}\cong \G_m^2\sqcup \mathbb{A}^1\sqcup \mathbb{A}^1
\end{equation}
where the three corresponding rulings are drawn in \Cref{fig:Front_31_ruling}. The identity \eqref{eq:ruling_stratification} is understood as a finite partition into locally closed subvarieties, cf.~\cite[Theorem 3.4]{HenryRutherford15}. The maximal normal ruling is shown in \Cref{fig:Front_31_ruling} (left) and determines the ruling piece $\G_m^2\cong\{(x,y)\in\mathbb{A}^2:x(1+xy)\neq0\}$. The submaximal ruling shown in \Cref{fig:Front_31_ruling} (center) determines the ruling piece $\mathbb{A}^1\cong\{(0,y,1):y\in\mathbb{A}^1\}$, and \Cref{fig:Front_31_ruling} (right) yields the affine line $D_z\cong \mathbb{A}^1\cong\{(1,-1,z):z\in\mathbb{A}^1\}$ as in \eqref{eq:affine_zline_trefoil}. 

\noindent (2) In these coordinates $(x,y,z)\in\mathbb{A}^3$, the proof of \Cref{lem:sheaves_trefoil} implies that a sheaf $\SF\in\M_1(\La_{3_1})$ has acyclic stalk $C$ in the central region if and only if $1+xy\neq0$. This locus $\{1+xy\neq0\}\times B\bG_m\sse \M_1(\La_{3_1})$ is the complement $\M_1(\La_{3_1})\setminus (D_z\times B\G_m)$ of the affine line
\begin{equation}\label{eq:affine_zline_trefoil}
D_z\times B\bG_m:=\{x=1, y=-1\}\times B\bG_m\sse \M_1(\La_{3_1})
\end{equation}
which is parametrized by the coordinate $z\in\mathbb{A}^1$. Note that $\M_1(\La_{3_1})\setminus D_z$ is isomorphic to
$$\{(x,y)\in\mathbb{A}^2: 1+xy\neq0\}\times B\bG_m.$$
From the perspective of cluster varieties, $1+xy$ is one of the cluster variables for the $A_2$-cluster algebra given by the ring of functions of $\{(x,y,z)\in\mathbb{A}^3: xyz+x+z=1\}$, and localizing at it leads to the $A_1$-cluster algebra with one frozen vertex, which is precisely $\{(x,y)\in\mathbb{A}^2: 1+xy\neq0\}$.

\subsection{Ruling decomposition for the moduli stack of the convex end}\label{ssec:ruling_decomposition_convex} 
\begin{center}
	\begin{figure}[H]
		\centering
		\includegraphics[scale=0.75]{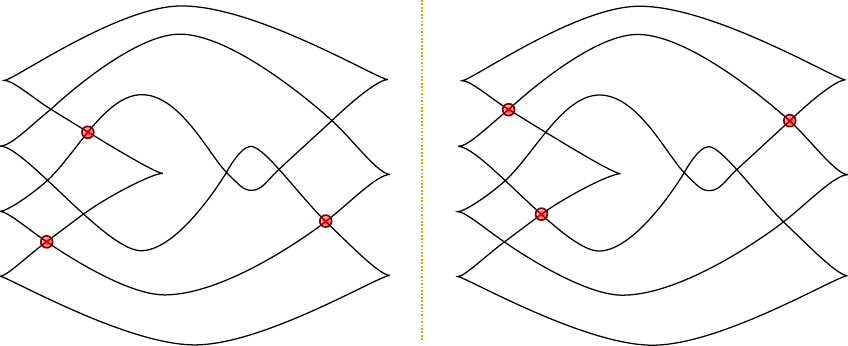}
		\caption{The two $\Z$-graded normal rulings for $\La_{m(9_{46})}$ corresponding to two embedded Lagrangian 2-disks bounding $\La_{m(9_{46})}$, cf.~\cite[Example 2.8]{CornwellNgSivek16}.}
		\label{fig:m946_front_rulings}
	\end{figure}
\end{center}

In line with \Cref{lem:sheaves_trefoil}, the purpose of this subsection is to describe the derived stack $\M(W(\La_{m(9_{46})}))$ of microlocal rank 1 sheaves with singular support in the Legendrian Whitehead double $W(\La_{m(9_{46})})$ of the unique max-tb Legendrian representative $\La_{m(9_{46})}$ of the smooth knot $m(9_{46})$. Fronts for $\La_{m(9_{46})}$ and $W(\La_{m(9_{46})})$ are respectively depicted in \Cref{fig:Front_SatelliteConcordance1} (right) and \Cref{fig:Front_SatelliteConcordance2} (right). Recall that in this article we use the word ruling to always mean a $\Z$-graded normal ruling. The two rulings for the Lagrangian 2-disks bounding $\La_{m(9_{46})}$ are depicted in \Cref{fig:m946_front_rulings}, see also \cite[Figure 3]{CornwellNgSivek16}. To ease notation we write $\M:=\M_1(W(\La_{m(9_{46})}))$. Rather than obtaining a complete description of $\M$, which would require a more significant amount of computation, we focus on studying its ruling decomposition.

\begin{center}
	\begin{figure}[H]
		\centering
		\includegraphics[width=\textwidth]{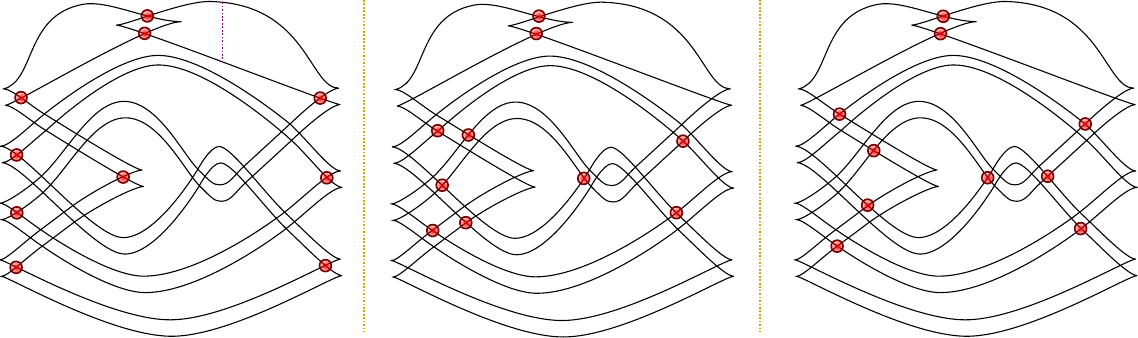}
		\caption{The three maximal $\Z$-graded normal rulings for $W(\La_{m(9_{46})})$. Each of the ten switches of each ruling has been marked with a red circle around it.}
		\label{fig:Wm946_initial_rulings}
	\end{figure}
\end{center}

\begin{center}
	\begin{figure}[H]
		\centering
		\includegraphics[width=\textwidth]{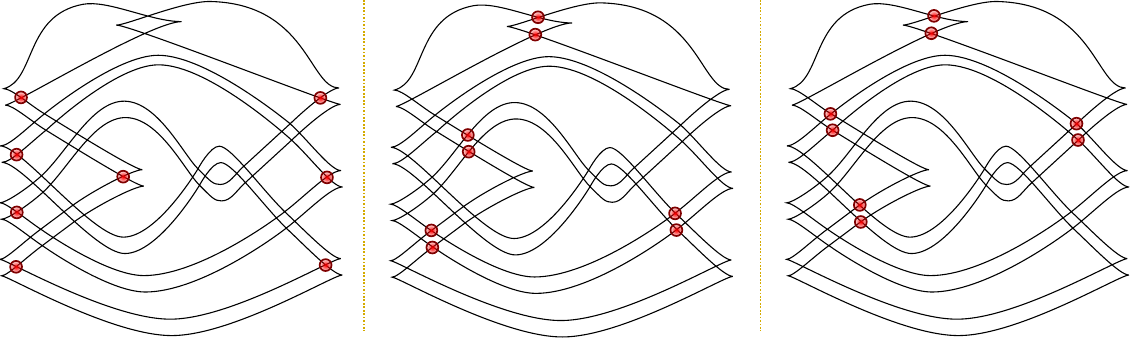}
		\caption{All the non-maximal rulings that we show connect to the maximal normal ruling in \Cref{fig:Wm946_initial_rulings} (left). Each of them has eight switches, marked with a red circle around each switch. Note that the rulings in the middle and right are exchanged by (what remains of) the $\Z_2$-symmetry of \Cref{fig:Front_SatelliteConcordance1} (right).}
		\label{fig:Wm946_initial_rulings2}
	\end{figure}
\end{center}

\begin{center}
	\begin{figure}[H]
		\centering
		\includegraphics[scale=1]{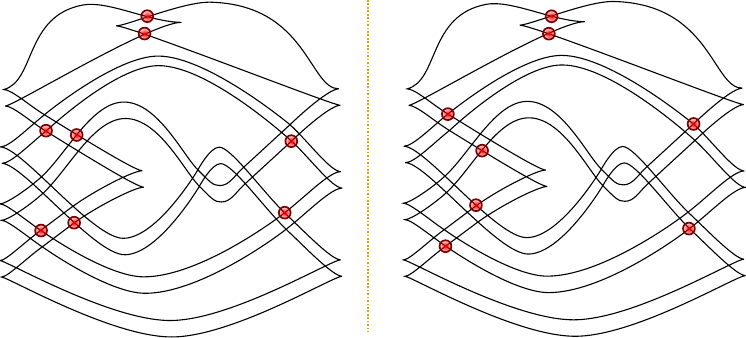}
		\caption{(Left) The unique non-maximal normal ruling $I_1\cong\bA^1$ that is a normal subruling of \Cref{fig:Wm946_initial_rulings} (center). (Right) The unique non-maximal normal ruling $I_2\cong\bA^1$ that is a normal subruling of \Cref{fig:Wm946_initial_rulings} (right).}
		\label{fig:Wm946_initial_rulings3}
	\end{figure}
\end{center}

\noindent The ruling decomposition of $\M$, using \Cref{fig:Front_SatelliteConcordance2} (right), consists of 8 pieces. Specifically, there are 3 maximal rulings, as depicted in \Cref{fig:Wm946_initial_rulings}, and eight rulings in total, depicted in \Cref{fig:Wm946_initial_rulings,fig:Wm946_initial_rulings2,fig:Wm946_initial_rulings3}. Here is our notation:
\begin{enumerate}
    \item The three maximal rulings in \Cref{fig:Wm946_initial_rulings} are denoted by $U^e,U^i_1,U^i_2$ left to right. At the 0-truncation level, each of them provides a homonymous piece of the ruling decomposition isomorphic to $\G_m^2$. The superscripts $e$ and $i$ are meant to indicate embedded and immersed.\\

    \item The three rulings in \Cref{fig:Wm946_initial_rulings2} are denoted by $L,D_1,D_2$ from left to right. At the 0-truncation level, each of them provides a a stratum (denoted by the same symbol) of the ruling decomposition which is isomorphic to $\bA^1$. Note that $D_1,D_2$ are obtained from the original rulings in \Cref{fig:m946_front_rulings} by doubling the corresponding switches and adding the two switches at the clasp. The ruling corresponding to $L$ is obtained from the standard ruling of the Reeb 2-copy, having switches at all crossings near doubled cusps, and extending it without any switches at the Whitehead clasp.\\

    \item The two rulings in \Cref{fig:Wm946_initial_rulings3} are denoted by $I_1,I_2$ left to right. At the 0-truncation level, each of them provides a homonymous piece of the ruling decomposition isomorphic to $\bA^1$.
\end{enumerate}

As a consistency check for the rulings above, the $\Z$-graded ruling polynomial of $W(\La_{m(9_{46})})$ is
\begin{equation}\label{eq:ruling_poly}
R_{W(\La_{m(9_{46})})}(z)=5z^{-1}+3z
\end{equation}
under the conventions of \cite{HenryRutherford15}. Indeed, the polynomial \eqref{eq:ruling_poly} is readily computed from the front \Cref{fig:Front_SatelliteConcordance2} (right) by applying the skein relation to the Whitehead clasp, which leads to two feasible cases.\footnote{One is the max-tb Legendrian unknot, which contributes $z^{-1}$ to \eqref{eq:ruling_poly}, the other the Reeb 2-copy of $\La_{m(9_{46})}$.} The rulings contributing to $z^{-1}$ each have eight switches, and there are 5 of them as illustrated in \Cref{fig:Wm946_initial_rulings2} and \Cref{fig:Wm946_initial_rulings3}, and the (maximal) rulings contributing to $z$ each have ten switches, and there are 3 of them, as drawn in \Cref{fig:Wm946_initial_rulings}.


\subsection{Ruling images under sheaf quantization} Let $C:\La_-\prec\La_+$ be an exact Lagrangian concordance and $W(C):W(\La_-)\prec W(\La_+)$ its Legendrian Whitehead double. Consider a sheaf $\SF_-\in\M_1(\La_-)$. By definition:

\begin{enumerate}

    \item We denote by $W_{sw}(\SF_-)\in \M_1(W(\La_-))$ the sheaf obtained as the cone of the zero multiple of the continuation map $c:T_{-\varepsilon}\SF_-\to T_{\varepsilon}\SF_-$, by extending the sheaf $\mbox{cone}(0\cdot c)$ (which again has singular support in $\La_-^{(2)}$) to the Whitehead clasp with two switches, one at each crossings of the clasp. For instance, in our cases of interest, the ruling associated to such a sheaf $W_{sw}(\SF_-)$ is as depicted in \Cref{fig:Front_31_ruling} (center).\\

    \item For $\lambda\in k^\times$, we denote by $W_{\lambda}(\SF_-)\in \M_1(W(\La_-))$ the sheaf obtained as the cone of $\lambda\cdot c$, where $c$ is the continuation map $c:T_{-\varepsilon}\SF_-\to T_{\varepsilon}\SF_-$, by extending the resulting sheaf $\mbox{cone}(\lambda\cdot c)$, which has singular support in $\La_-^{(2)}$, trivially above the Whitehead clasp, i.e.~with no switches at the two crossings of the Whitehead clasp. A ruling associated to such a sheaf $W_{\lambda}(\SF_-)$ is as depicted in \Cref{fig:Front_31_ruling} (right).
\end{enumerate}

The images of such sheaves under the sheaf quantization of $W(C)$ have a simple description:
\begin{lemma}\label{lem:sq}
Let $C:\La_-\prec\La_+$ be an exact Lagrangian concordance and $W(C):W(\La_-)\prec W(\La_+)$ its Legendrian Whitehead double. Consider a sheaf $\SF_-\in\M_1(\La_-)$. Then the following equivalences hold:
\begin{equation}\label{eq:sheaf_quantization_rulings}
\Phi_{W(C)}(W_{sw}(\SF_-))\cong W_{sw}(\Phi_C(\SF_-))\qquad\&\qquad \Phi_{W(C)}(W_{\lambda}(\SF_-))\cong W_{\lambda}(\Phi_C(\SF_-)).\qquad 
\end{equation}
\end{lemma}

\begin{proof}
For a Lagrangian concordance $L:\La_-\prec\La_+$ and $\SF_-\in\M_1(\La_-)$, the image $\Phi_L(\SF_-)$ is given by the restriction to the convex end of the unique sheaf with singular support on the Legendrian lift of $L$ which restricts to $\SF_-$ at its concave end (and has the corresponding local system specified there). See \cite[Theorems 1.2\&1.3]{li2024lagrangiancobordismfunctormicrolocal}. In our case, consider the Legendrian Whitehead doubling as two steps: first taking the Reeb 2-copy and then adding the Whitehead clasp.

\noindent For the first step, suppose that $\SF_C$ is a sheaf quantization of the concordance $C$. Then the relevant sheaf quantizations of the Reeb 2-copy $C^{(2)}$ of $C$ are given by cones of scalar multiples of the continuation map $T_{-\varepsilon}\SF_C\to T_{\varepsilon}\SF_C$, cf.~\cite[Prop.~3.10]{li2024lagrangiancobordismfunctormicrolocal}. If it is the zero multiple, then we obtain the direct sum sheaves on both ends, and if it is non-zero, then we obtain extension sheaves on both ends. In the former case, the ruling associated to the sheaf on the convex end of $C^{(2)}$ is obtained by doubling the switches of the ruling for $\Phi_C(\SF_-)$. In the latter case, the rulings are those having switches exactly at the crossings near the doubled cusps. See \Cref{fig:Wm946_initial_rulings_images}. For the second step, adding the Whitehead clasp, we can sheaf quantize around the tangle along the satellited concordance by trivially extending the sheaf quantization on the concave end throughout to the convex end. This second step leaves the sheaf and its associated ruling identical: if it had no switches (resp.~two switches) at the Whitehead clasp in the concave end, then it will have no switches (resp.~two switches) in the convex end. Putting both steps together concludes the statement.
\end{proof}


\subsection{Non-decomposability of the Cornwell-Ng-Sivek concordance}\label{ssec:non_decomposability_CNS} The goal of this subsection is to prove \Cref{thm:main}. Let $C_\pm:\La_\st\prec\La_{m(9_{46})}$ be the Lagrangian concordances obtained by puncturing the embedded Lagrangian 2-disk bounding $\La_{m(9_{46})}$ with $\Z$-graded normal ruling as depicted in \Cref{fig:m946_front_rulings} (left), for $C_+$, and \Cref{fig:m946_front_rulings} (right) for $C_-$. Let $C^w_\pm:\La_{3_1}\prec W(\La_{m(9_{46})})$ denote the Cornwell-Ng-Sivek concordances, obtained by satelliting $C_{\pm}$ along the Legendrian Whitehead double pattern as in \cite[Figure 8]{CornwellNgSivek16}. Their ends are depicted in \Cref{fig:Front_SatelliteConcordance2}. Denote by
\begin{equation}\label{eq:induced_map_CNS}
\Phi_{C^w_\pm}:\M_1(\La_{3_1})\lr \M_1(W(\La_{m(9_{46})}))
\end{equation}
their induced maps under sheaf quantization, cf.~\cite[Theorem 1.1]{li2024lagrangiancobordismfunctormicrolocal}, which are always open immersions. In order to prove \Cref{thm:main}, we proceed by contradiction. For that, we assume that either of $C^w_\pm$ is decomposable, say $C^w_+$ for specificity. By \cref{lem:decomposable_closed},
\begin{equation}\label{eq:induced_map_CNS_plus}
\Phi_{C^w_+}:\M_1(\La_{3_1})\lr \M_1(W(\La_{m(9_{46})}))
\end{equation}
must then be an open and closed immersion.\\

Let us denote the three rulings for $\M_1(\La_{3_1})$ in \Cref{fig:Front_31_ruling} by $T,A_1,A_2$, reading left to right. The two submaximal rulings $A_1,A_2$, whose pieces are isomorphic to $\bA^1$, are depicted in the left column of \Cref{fig:Wm946_initial_rulings3} as well. The images of these ruling pieces under $C^w_\pm$ are depicted in \Cref{fig:Wm946_initial_rulings_images}, according to \Cref{lem:sq}. The figure is to be interpreted as follows: for any sheaf $\SF\in \M_1(\La_{3_1})$ with ruling given by the first column, the sheaf $\Phi_{C^w_\pm}(\SF)$ has associated ruling given by the second or third columns.

\begin{center}
	\begin{figure}[H]
		\centering
		\includegraphics[width=\textwidth]{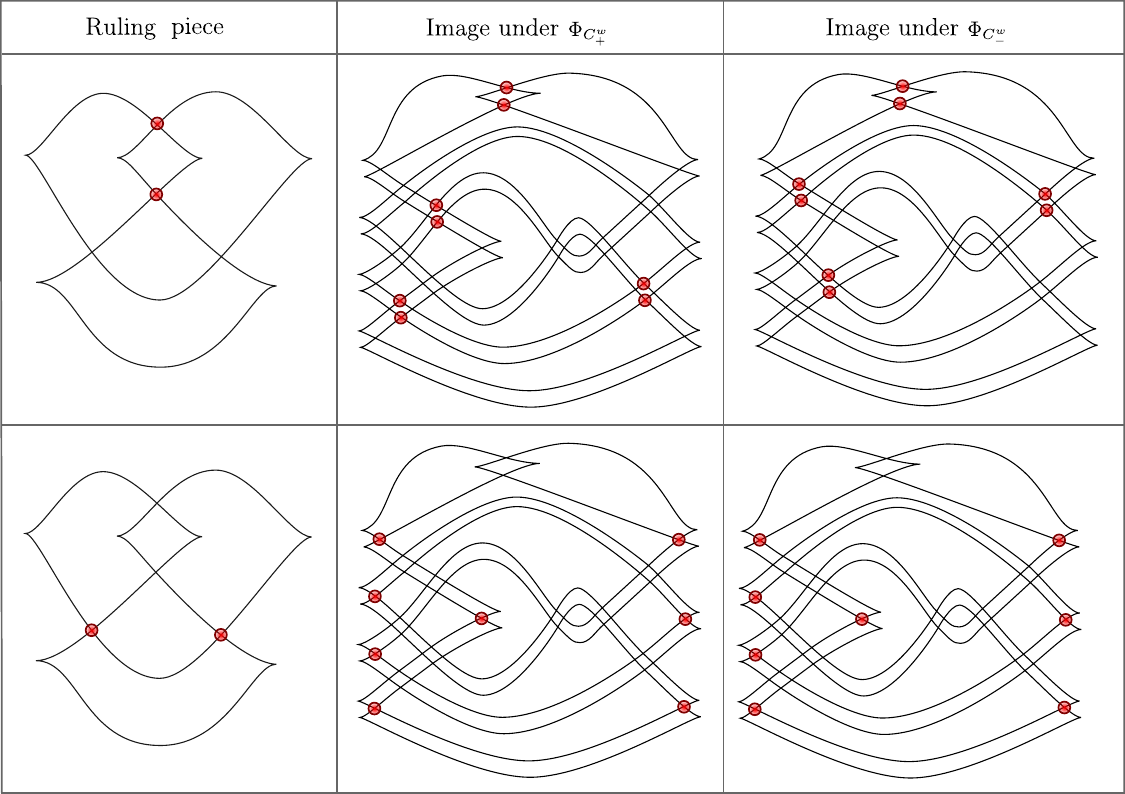}
		\caption{(Left) The two non-maximal ruling strata for this front of the max-tb Legendrian trefoil knot. (Center) The corresponding images of the rulings in the left column under the sheaf quantization $\Phi_{C^w_+}$. (Right) The images of the rulings in the left column under the other sheaf quantization $\Phi_{C^w_-}$.}
		\label{fig:Wm946_initial_rulings_images}
	\end{figure}
\end{center}

In order to show that \eqref{eq:induced_map_CNS_plus} cannot be an open and closed immersion, we use an Euler characteristic count. Specifically, we consider an algebraically closed field $k$ of characteristic zero and let $\chi_c$ be the compactly supported Euler characteristic for stacks introduced in \cite[Definition 4.8]{Joyce_EulerCharacteristic}. By \cite[Lemma 4.13]{Joyce_EulerCharacteristic} it is additive under disjoint unions and multiplicative under Cartesian products. We will use\footnote{In general a piece $\bA^n\times\G^{g}_m\times BG$ of a finite ruling decomposition contributes either one, if $g=0$, or zero otherwise.} the three values
\begin{equation}\label{eq:values_chi}
\chi_c(\bA_1)=1,\qquad\chi_c(\G_m)=0,\qquad \chi_c(B\bG_m)=1.
\end{equation}
In particular, the ruling decomposition gives
\begin{equation}\label{eq:values_chi_concave}
\chi_c(\M_1(\La_{3_1}))=\chi_c(\G_m^2\times B\bG_m)+2\chi_c(\bA^1\times B\bG_m)=\chi_c(\G_m^2)+2\chi_c(\bA^1)=0+2=2.
\end{equation}
\noindent In the ruling decomposition of the target of \eqref{eq:induced_map_CNS_plus}, all ruling pieces contribute nonnegatively:
\begin{equation}\label{eq:values_non_negative}
\chi_c(\S)\geq0,\quad \mbox{ if }\S\sse \M_1(W(\La_{m(9_{46})}))\mbox{ is a ruling piece}.
\end{equation}
\noindent Let $\Im\sse \M_1(W(\La_{m(9_{46})}))$ denote the image of \eqref{eq:induced_map_CNS_plus}, which is an open and closed subset. Now, by \Cref{lem:sq}, the ruling piece $A_1\sse \M_1(\La_{3_1})$ maps to the ruling piece $D_1$, cf.~\Cref{ssec:ruling_decomposition_convex}, and the ruling piece $A_2\sse \M_1(\La_{3_1})$ maps to the ruling piece $L$, see also \Cref{fig:Wm946_initial_rulings_images}. Since the intersection $\Im\cap L$ is open and closed and $L$ is connected, we conclude that $L\sse\Im$. Similarly, the intersection $\Im\cap D_1$ is open and closed and $D_1$ is connected, so we also have the inclusion $D_1\sse\Im$. We record these inclusions:
\begin{equation}\label{eq:inclusions1}
L\sse\Im\qquad\&\qquad  D_1\sse\Im
\end{equation}
\noindent Let us now argue that we must also have $D_2\sse\Im$. Indeed, since $\Phi_{C^w_-}$ is continuous, the pre-image $\Phi^{-1}_{C^w_-}(\Im)\sse \M_1(\La_{3_1})$ is an open and closed subset. This pre-image is not empty because $L\sse\Im$, as in \eqref{eq:inclusions1}, and $\Phi_{C^w_-}$ (also) maps the ruling piece $A_2$ to $L$, by \Cref{lem:sq}, cf.~\Cref{fig:Wm946_initial_rulings_images}. Since $\M_1(\La_{3_1})$ is connected, it must be that $\Phi^{-1}_{C^w_-}(\Im)=\M_1(\La_{3_1})$. Given that $\Phi_{C^w_-}$ maps the ruling piece $A_1$ to the ruling piece $D_2$, $D_2$ intersects the image of $\Phi_{C^w_-}$, cf.~\Cref{fig:Wm946_initial_rulings_images}. As before, $D_2$ being connected and the intersection $D_2\cap\Im$ being open and closed leads to the inclusion
\begin{equation}\label{eq:inclusions2}
D_2\sse\Im
\end{equation}
The image $\Im\sse \M_1(W(\La_{m(9_{46})}))$ being open and closed implies that it either entirely contains a piece of the ruling decomposition of $\M_1(W(\La_{m(9_{46})}))$ or it is disjoint from it. By using \eqref{eq:inclusions1} and \eqref{eq:inclusions2} we can write
\begin{equation}\label{eq:inclusions3}
\Im= L\sqcup D_1\sqcup D_2\sqcup \S
\end{equation}
where $\S$ is a (possibly empty) disjoint union of some of the remaining pieces of the ruling decomposition of $\M_1(W(\La_{m(9_{46})}))$. By \eqref{eq:values_non_negative}, $\chi_c(\S)\geq0$, i.e.~it must be non-negative, and therefore \eqref{eq:inclusions3} implies
\begin{equation}\label{eq:values_chi_image}
\chi_c(\Im)=\chi_c(L)+\chi_c(D_1)+\chi_c(D_2)+\chi_c(\S)\geq \chi_c(L)+\chi_c(D_1)+\chi_c(D_2)=3.
\end{equation}
Since $\Im$ is the image of an open and closed immersion, $\chi_c(\Im)=\chi_c(\M_1(\La_{3_1}))$, which leads to the contradiction
$$3\stackrel{\eqref{eq:values_chi_image}}{\leq}\chi_c(\Im)=\chi_c(\M_1(\La_{3_1}))\stackrel{\eqref{eq:values_chi_concave}}{=}2.$$
\noindent This concludes that the open immersion $\Phi_{C^w_+}$ cannot also be a closed immersion. By \Cref{lem:decomposable_closed}, $C^w_+$ cannot be a decomposable Lagrangian concordance. Finally, the same argument, interchanging the roles of $C^w_\pm$, concludes that $C^w_-$ is not decomposable either.\hfill$\Box$


\subsection{Non-decomposability of the Reeb 2-copy}\label{ssec:non_decomposability_Reeb2copy} The goal of this subsection is to prove \Cref{thm:main2}. We proceed by contradiction, and thus we assume that the Reeb 2-copy $C^{(2)}_+:\La_H\prec \La_{m(9_{46})}^{(2)}$ is decomposable. By \Cref{lem:decomposable_closed}, its sheaf quantization
\begin{equation}\label{eq:Reeb_2_copy_1}
\Phi_{C^{(2)}_+}:\M_1(\La_H)\lr\M_1(\La_{m(9_{46})}^{(2)})
\end{equation}
must be an open and closed immersion. In this case we will exploit the $\Z_2$-symmetry of the front of $\La_{m(9_{46})}^{(2)}$ inherited from that of the front of $\La_{m(9_{46})}$, i.e.~the symmetry under the map $(x,z)\mapsto (x,-z)$.

\begin{center}
	\begin{figure}[H]
		\centering
		\includegraphics[width=\textwidth]{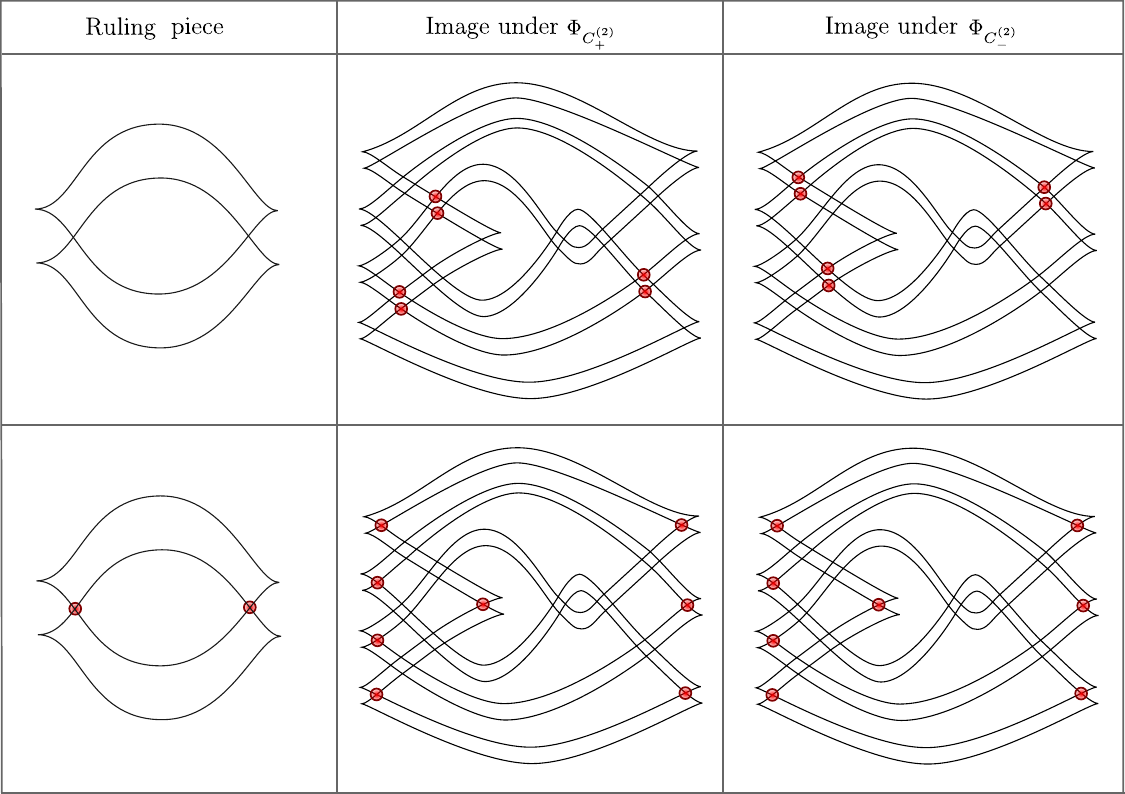}
		\caption{(Left) The two ruling strata for this front of the standard Legendrian Hopf link $\La_H$, realized as the Reeb 2-copy of the max-tb Legendrian unknot. (Center) The corresponding images of the rulings in the left column under the sheaf quantization $\Phi_{C^{(2)}_+}$. (Right) The images of the rulings in the left column under the other sheaf quantization $\Phi_{C^{(2)}_-}$.}
		\label{fig:Wm946_initial_rulings_images}
	\end{figure}
\end{center}

For that, let $I_+\sse \M(\La_{m(9_{46})}^{(2)})$ be the ruling piece associated to the ruling in the center of the first row of \Cref{fig:Wm946_initial_rulings_images}, which is obtained by doubling the switches of the ruling for the corresponding Lagrangian 2-disk, as drawn in \Cref{fig:m946_front_rulings} (left). Similarly, we denote by $I_-$ and $P$ the ruling pieces associated to the rulings at the right of the first row of \Cref{fig:Wm946_initial_rulings_images} and at the right (or center) of the second row. By construction, the $\Z_2$-symmetry induces an involution $\eta:\M(\La_{m(9_{46})}^{(2)})\lr\M(\La_{m(9_{46})}^{(2)})$ which exchanges $I_\pm$, i.e.~$\eta(I_+)=I_-$, and preserves $P$, i.e.~$\eta(P)=P$, see \Cref{fig:Wm946_initial_rulings_images}. In particular, the image $\Im:=\im\Phi_{C^{(2)}_+}$ is an open and closed subset which intersects the $\Z_2$-invariant subset $P$, cf.~\Cref{lem:sq}. Since both $\Im$ and $\eta(\Im)$ are open and closed, connected, and their intersection is non-empty, as it contains $P$, we conclude that $\Im=\eta(\Im)$. Given that we are assuming that $\Phi_{C^{(2)}_+}$ in \eqref{eq:Reeb_2_copy_1} is open and closed, it must be an isomorphism onto its image. In particular, the $\eta$-invariance of its image provides a $\Z_2$-symmetry of its domain $\M_1(\La_H)$, i.e.~an involution $\eta':\M_1(\La_H)\lr \M_1(\La_H)$.

\noindent By \cite[Section 2.5.2]{casals2024positivemicrolocalholonomiesglobally}, the derived stack $\M_1(\La_H)$ can be expressed as the global quotient
$$\M_1(\La_H)\cong \left[\{(x,y)\in\bA^2: 1+xy\neq0\}/\G_m^2\right],\mbox{ with }\G_m^2\mbox{-action }(t_1,t_2)\cdot (x,y)=(t_1t_2^{-1}x,t_2t_1^{-1}y).$$
\noindent In this description, the only point in $\M_1(\La_H)$ with $\G_m^2$-stabilizer is the origin $O:=(0,0)/\G_m^2$. Therefore it must be fixed under the involution $\eta'$, i.e.~$\eta'(O)=O$. Now, the image of the origin under \eqref{eq:Reeb_2_copy_1} lies in the ruling stratum $I_+$, cf.~\Cref{lem:sq} and first row of \Cref{fig:Wm946_initial_rulings_images}. This is a contradiction, because \eqref{eq:Reeb_2_copy_1} is $(\eta',\eta)$-equivariant by construction, $\eta(I_+)=I_-$ and $I_+$ is disjoint from $I_-$. Indeed, $(\eta',\eta)$-equivariance would imply
$$I_-=\eta(I_+)=\eta(\Phi_{C^{(2)}_+}(O))=\Phi_{C^{(2)}_+}(\eta'(O))=\Phi_{C^{(2)}_+}(O)=I_+,$$
while the ruling decomposition gives $I_+\cap I_-=\emptyset$. We thus conclude that the open immersion $\Phi_{C^{(2)}_+}$ cannot be an open and closed immersion. By \Cref{lem:decomposable_closed}, the concordance $C^{(2)}_+$ is not decomposable. A symmetric argument shows that $C^{(2)}_-$ is not decomposable either. \hfill$\Box$

\begin{remark}
In \Cref{ssec:non_decomposability_CNS} we obstructed $\Phi_{C^{w}_+}$ from being an open and closed immersion by using Euler characteristics. In \Cref{ssec:non_decomposability_Reeb2copy} we obstructed $\Phi_{C^{(2)}_+}$ by using some information about the stabilizers and the $\Z_2$-symmetry. In general, one can obstruct a morphism of derived stacks from being an open and closed immersion in a variety of ways. For instance, one can try to compare the $E$-polynomials (Hodge-Deligne polynomials) of the domain with those of the image, explore base changes, or use finer cohomological computations. In general, being an open and closed immersion forces an isomorphism onto an open and closed substack: among all the properties to obstruct, this one is in principle relatively feasible if one has some information on the domain and its codomain.\sbs
\end{remark}


\subsection{Final comments}\label{ssec:final_comments} The arguments of \Cref{ssec:non_decomposability_CNS} and \Cref{ssec:non_decomposability_Reeb2copy} do not require a complete understanding of the isomorphism type of the sheaf moduli $\M_1$ of their convex ends. It is nevertheless possible to describe them by studying the possible sheaves on the front. For $\M:=\M_1(W(\La_{m(9_{46})}))$ the description is as follows. The derived stack $\M$ has three connected components, each irreducible. The first component consists of
\begin{equation}\label{eq:iso_type_CNS_concave}
\overline{U^e}\times B\bG_m\cong \left(U^e\sqcup L\sqcup D_1\sqcup D_2\right)\times B\bG_m\cong \left(\SA_2\cup_{\SA_1}\SA_2\right)\times B\bG_m.
\end{equation}
where $U^e$ is the stratum associated to the leftmost maximal ruling in \Cref{fig:Wm946_initial_rulings} and the strata for $L,D_1,D_2$ are associated to the three rulings in \Cref{fig:Wm946_initial_rulings2}. In \eqref{eq:iso_type_CNS_concave} we have denoted
\begin{equation}\label{eq:notation_varieties}
\SA_1:=\{(x,y)\in\mathbb{A}^2_k: xy+1\neq0\},\quad \SA_2:=\{(x,y,z)\in\mathbb{A}^3_k: xyz+x+z=1\}
\end{equation}
considered as derived stacks, by viewing an affine scheme as a (rather trivial) derived stack.\footnote{The notation is motivated by the fact that they carry $A_1$ and $A_2$ cluster structures, respectively.} Note that we have a natural open embedding
\begin{equation}\label{eq:inclusionA1A2}
\iota:\SA_1\hookrightarrow \SA_2,\qquad (x,y)\longmapsto \left(x,y,\frac{1-x}{1+xy}\right).
\end{equation}
\noindent In \eqref{eq:iso_type_CNS_concave} we denoted by $\SA_2\cup_{\SA_1}\SA_2$ the non-separated scheme given by gluing two copies of $\SA_2$ along the open set $\SA_1$. That is, $\SA_2\cup_{\SA_1}\SA_2$ is the homotopy pushout of $\iota$ with a copy of itself. Thus $\SA_2\cup_{\SA_1}\SA_2$ is isomorphic to $\SA_2$ except that the line $D_z\sse\SA_2$ has been doubled. In summary, the component indexed by the maximal ruling $U^e$, where all the embedded Lagrangian fillings lie, is isomorphic to the derived stack of the trefoil but with a doubled line. Namely, there are two copies of $\SA_2$, realized with the locally closed partitions
$$\SA_2\cong U^e\sqcup L\sqcup D_1\cong U^e\sqcup L\sqcup D_2$$
and they are glued along the complements of the lines $D_1$ and $D_2$ respectively. Specifically, the lines $D_1,D_2$ act as a double line on the complement of $\SA_1\cong U^e\sqcup L$.

\noindent The connected component $\overline{U^e}$ is the one where the image of any embedded Lagrangian filling must lie. There are nevertheless two more connected components, which keep track of unobstructed immersed Lagrangian fillings. Specifically, the second and third connected components of $\M$ consist of
$$\overline{U^i_1}\cong U^i_1\sqcup I_1\cong \SA_1\times B(\G_m\times \G_a),\qquad \overline{U^i_2}\cong U^i_2\sqcup I_2\cong\SA_1\times B(\G_m\times \G_a),$$
where $U_1^i,U_2^i$ denote the strata associated to the maximal rulings in \Cref{fig:Wm946_initial_rulings} (middle and right). Note that the existence of the stabilizer $\G_m\times\G_a$ implies that an open and closed immersion from $\M_1(\La_{3_1})$ must necessarily land in the connected component $\overline{U^e}\times B\G_m$ in \eqref{eq:iso_type_CNS_concave}. Therefore, by the Ax-Grothendieck theorem, cf.~\cite[Theorem]{Ax69}, there can be no open and closed immersion from $\A_2\times B\G_m$, which is the microlocal rank one moduli for the concave end, to this connected component $\overline{U^e}\times B\G_m$. In particular, we obtain the following strengthening of \Cref{thm:main}:

\begin{cor}\label{eq:cor_not_any}
There exists no decomposable Lagrangian concordance from $\La_{3_1}$ to $W(\La_{m(9_{46})})$.
\end{cor}

A similar description can be given for $\M_1(\La^{(2)}_{m(9_{46})})$ and the analogous strengthening of  \Cref{thm:main2} also holds. Specifically, this moduli stack also has three connected components, each irreducible. For instance, two of the components of the stack (Postnikov) 1-truncation of $\M_1(\La^{(2)}_{m(9_{46})})$ are isomorphic to the quotient stack
\begin{equation}\label{eq:component_doubled}
[\SA_1/B_-],\mbox{ where }B_-\mbox{ acts by }
\begin{pmatrix}
a & 0\\
b & c
\end{pmatrix}\begin{pmatrix}
x\\
y
\end{pmatrix}=\begin{pmatrix}
cx/a\\
ay/c
\end{pmatrix}
\end{equation}
and $B_-$ denotes the Borel subgroup of lower triangular matrices.

\begin{remark}
If we consider the untruncated derived stack $\M_1(\La^{(2)}_{m(9_{46})})$, then there should be additional contributions to $\pi_2$, coming from the fact that some of the sheaves in these two components are objects with non-trivial Ext-groups in degree -1.\sbs
\end{remark}

\noindent The third component, similar to \eqref{eq:iso_type_CNS_concave} above, is isomorphic to the sheaf moduli on the concave end with a doubled divisor. Specifically, it is isomorphic to
\begin{equation}\label{eq:component_doubled2}
\left[\SA_1\cup_{\{x\neq0\}}\SA_1\right]/\G_m^2,
\end{equation}
with the diagonal torus action as in \eqref{eq:component_doubled} above. Then the analogous statement of \Cref{eq:cor_not_any} holds for \Cref{thm:main2}. It is my understanding that, after this manuscript was shared with them, J.~Ma, L.~Ng, J.~Hughes and Y.~Pan are working on such results and their improvements, from a Floer-theoretic viewpoint. We therefore do not delve further into the details of these extensions, as they will be part of their upcoming work.\\

Finally, the following comments stem from conversations with L.~Ng and S.~Sivek, who kindly shared their thoughts after I sent them a first draft of the manuscript:

\begin{enumerate}[label=$(\arabic*)$]
    \item Suppose that a Lagrangian concordance $C:\La_1\prec\La_2$ is such that its sheaf quantization $\Phi_C$ is {\it not} a closed immersion, and thus $C$ is not decomposable. Then for any decomposable Lagrangian concordance $L:\La_2\prec\La_3$, the concatenation $L\circ C:\La_1\prec\La_3$ is not decomposable either. Indeed, since $\Phi_L$ is an open and closed immersion, the composition $\Phi_{L\circ C}=\Phi_L\circ\Phi_C$ is an open and closed immersion if and only if $\Phi_C$ is an open and closed immersion.\\

    \item The concatenation of a decomposable Lagrangian cobordism {\it below} a not necessarily decomposable Lagrangian concordance can result in a decomposable cobordism, in partial contrast to $(1)$ above. For instance, consider the non-decomposable Lagrangian concordance from \Cref{thm:main2} and assume that we have performed the decomposable concordance from $\La_\st$ to $\La_{m(9_{46})}$ so that a point $p\in\La_\st$ stays fixed throughout it, which can always be achieved. Let $L:\emptyset\prec\La_H$ denote the decomposable Lagrangian filling obtained by pinching at a short vertical arc right above $p$, joining the two components of the Reeb 2-copy $\La_H$ of the max-tb unknot, and filling the resulting max-tb unknot with a minimum cobordism. Then the resulting Lagrangian filling
    \begin{equation}\label{eq:filling_Hopf}
        C\circ L:\emptyset\prec \La_{m(9_{46})}^{(2)}
    \end{equation}
    is decomposable, even if $C$ is not. It is decomposable because $C\circ L$ is Hamiltonian isotopic, relative to its convex end, to the Lagrangian filling of its convex end obtained by pinching at the short vertical arc above the image of $p$, joining the two components of the Reeb 2-copy of $\La_{m(9_{46})}$ and filling the resulting max-tb unknot with a minimum cobordism. Indeed, a Hamiltonian isotopy is given by the family of fronts obtained by sliding the pinching move at $p$ up and down the satellited concordance.
\end{enumerate}

\begin{remark}
$(i)$ A similar construction illustrates the phenomenon in $(2)$ for the non-decomposable Lagrangian in \Cref{thm:main}. Specifically, consider the Lagrangian filling $L:\emptyset\prec\La_{3_1}$ of the max-tb right-handed trefoil $\La_{3_1}$ obtained from the front in \Cref{fig:Front_SatelliteConcordance2} (left) as follows. First we pinch at a short vertical arc (close) to the right of the right cusp of the Whitehead clasp, and then we pinch at a diagonal arc from the top strand of that right cusp in the clasp to the bottom strand of the left cusp in the clasp. These two pinchings result in a max-tb unknot, and we take $L$ to be the Lagrangian filling obtained by further filling the unknot with a minimum cobordism. Similar to $(2)$ above, the composition of $L$ with the Lagrangian concordance $C$ from \Cref{thm:main} results in a Lagrangian filling
   \begin{equation}\label{eq:filling_trefoil}
        C\circ L:\emptyset\prec W(\La_{m(9_{46})})
    \end{equation}
that, like \eqref{eq:filling_Hopf}, is decomposable despite $C$ not being decomposable. Indeed, the Lagrangian filling given by \eqref{eq:filling_trefoil} is Hamiltonian isotopic, relative to its convex end, to the Lagrangian filling of $W(\La_{m(9_{46})})$ obtained by performing the same two pinchings near the Whitehead clasp and filling the resulting max-tb unknot with a minimum cobordism.\\

\noindent $(ii)$ The constructions in $(2)$ and $(i)$ above are rather general: as long as we have a decomposable Lagrangian filling $L_{pat}$ of a given pattern in $(J^1S^1,\xi_\st)$, the satellite of any Lagrangian concordance along that pattern provides two natural isotopic Lagrangian fillings of the satellited convex end. The first Lagrangian filling is obtained by filling the satellited concave end with $L_{pat}$ tand then concatenating with the satellited concordance. The second Lagrangian filling is obtained by directly using $L_{pat}$ at the satellited convex end of the satellited concordance. These two Lagrangian fillings are Hamiltonian isotopic relative to the convex end, with the Hamiltonian isotopy being given by sliding where $L_{pat}$ is implemented from one end of the satellited concordance to the other. Thus the former Lagrangian filling is always decomposable, even if its satellited Lagrangian concordance piece is not, because the latter Lagrangian filling is decomposable by construction.\sbs
\end{remark}

Finally, we conclude the article by emphasizing that it should be possible to construct non-decomposable Lagrangian fillings of the Reeb 2-copy or the Whitehead double of the max-tb $\La_{m(9_{46})}$ by filling the concave ends of these non-decomposable concordances with an appropriate chosen (even decomposable) Lagrangian fillings. This will be hopefully explained in detail in the upcoming work of J.~Ma, L.~Ng, J.~Hughes and Y.~Pan. In general, the message that we want to convey is that one should not expect the results of this article to be restricted simply to Lagrangian concordances, but rather the unraveling start of the fact that non-decomposable Lagrangian concordances, Lagrangian cobordisms, and Lagrangian fillings all abound.

\bibliography{main.bib}
\end{document}